\documentclass[12pt]{amsart}
\usepackage[all]{xy}

\swapnumbers
\newtheorem{thm}[subsection]{Theorem}
\newtheorem{prop}[subsection]{Proposition}
\newtheorem{lemma}[subsection]{Lemma}

\def\t{\otimes}

\def\r{\succ}
\def\l{\prec}
\def\g{\dashv}
\def\d{\vdash}

\def\Hom{\operatorname{Hom}}

\def\PP{{\mathcal P}}
\def\QQ{{\mathcal Q}}
\def\AA{{\mathcal A}}
\def\BB{{\mathcal B}}
\def\XX{{\mathcal X}}

\def\nw{\nwarrow}
\def\ne{\nearrow}

\def\sw{\swarrow}
\def\se{\searrow}

\def\sq{\ \square \ }

\begin{document}

\author[J.-L. Loday]{Jean-Louis Loday}
\address{Institut de Recherche Math\'ematique Avanc\'ee\\
    CNRS et Universit\'e Louis Pasteur\\
    7 rue R. Descartes\\
    67084 Strasbourg Cedex, France}
\email{loday@math.u-strasbg.fr}
\urladdr{www-irma.u-strasbg.fr/{$\sim$}loday/}
\title{Completing the operadic butterfly}

\date{\today}

\begin{abstract}
We complete a certain diagram (the operadic butterfly) of categories of algebras involving $Com$, $As$, and $Lie$ by constructing a type of algebras which have 4 generating operations and 16 relations. The associated operad is self-dual for Koszul duality.
\end{abstract}

\maketitle

\section{Introduction} \label{S:int} The three categories of algebras  $Com$, $As$, and $Lie$ (for commutative algebras, associative algebras and Lie algebras) are related by two functors
$$ Com \hookrightarrow As \stackrel{-}{\to} Lie \ .$$

There are other types of algebras with similar functors which make up the following ``operadic butterfly":

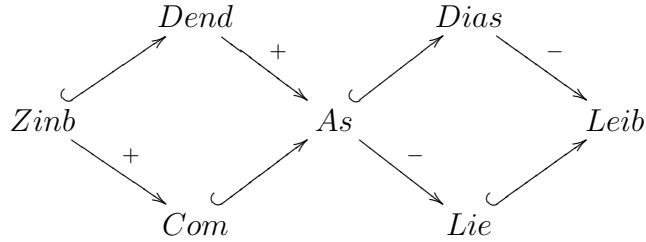
\[\xymatrix{
 &Dend\ar[rd]^{+} & &Dias\ar[rd]^{-} & \\
 Zinb\ar@{^{(}->}[ru]\ar[rd]^{+}& &As\ar@{^{(}->}[ru]\ar[rd]^{-} & &Leib \\
 &Com\ar@{^{(}->}[ru] & &Lie\ar@{^{(}->}[ru] & 
}\]
\\

The four other categories of algebras which appear in this diagram are as follows:

\noindent $Zinb = $ category of \emph{Zinbiel} algebras. They have one operation $x\cdot y$ (with no symmetry), satisfying 
$$(x\cdot y) \cdot z = x\cdot (y \cdot z) + x\cdot (z \cdot y).$$

\noindent $Dend = $ category of \emph{dendriform} algebras. They have two operations $x\l y$ and $x\r y$ (with no symmetry), satisfying 
\begin{displaymath}
\left\{\begin{array}{rcl}
(x\l y)\l z &=& x\l (y\l z) + x\l (y\r z)  , \\
(x\r y)\l z &=& x\r (y\l z) , \\
(x\l y)\r z  + (x\r y)\r z &=& x\r (y\r z).
\end{array}
\right.
\end{displaymath}

\noindent $Dias = $ category of \emph{diassociative} algebras (or associative dialgebras). They have  two operations $x\g y$ and $x\d y$ (with no symmetry), satisfying 
\begin{displaymath}
\left\{\begin{array}{c}
x\g (y\g z) = (x\g y)\g z = x\g (y\d z), \\
(x\d y)\g z = x\d (y\g z) , \\
(x\g y)\d z = x\d (y\d z) = (x\d y)\d z  . 
\end{array}
\right.
\end{displaymath}

\noindent $Leib = $ category of \emph{Leibniz} algebras. They have one operation $[x, y]$ (with no symmetry), satisfying $$[[x,y], z] = [[x, z], y] + [x,[y,z]].$$
\\

An arrow like $\AA \hookrightarrow\BB$ means that an $\AA$-algebra is a $\BB$-algebra satisfying some symmetry property. An arrow like $\AA\stackrel{\pm}{\to}\BB$ means that any generating operation of a $\BB$-algebra is obtained  by some addition (or substraction) of two operations of the $\AA$-algebra type (like the Lie bracket from the associative product). We let the reader to write the exact formulas in each case.

Each one of these types of algebras defines a binary quadratic operad. For these operads there is a well-defined notion of Koszul duality theory devised by Ginzburg and Kapranov \cite {G-K}. Let $\PP^!$ be the dual of the operad $\PP$ (note that $\PP^{!!} = \PP$). It turns out that Koszul duality in the operadic butterfly corresponds to symmetry around the vertical axis passing through $As$:
$$ As^! = As , \ Com ^! = Lie,\  Zinb^! = Leib, \ Dend^! = Dias \ . $$
A functor of the form $\hookrightarrow$ is changed into a functor of the form $\stackrel{\pm}{\leftarrow}$ by duality.

One can slightly enhance the operadic butterfly by putting the category $Vect$ of algebras with no operation (also called \emph{ abelian Lie algebras}) in between $Com$ and $Lie$. The operad $Vect$ is self-dual. So an immediate question comes to mind: can one complete the operadic butterfly by putting some category of algebras $\XX$ at the upper place on the middle axis ?

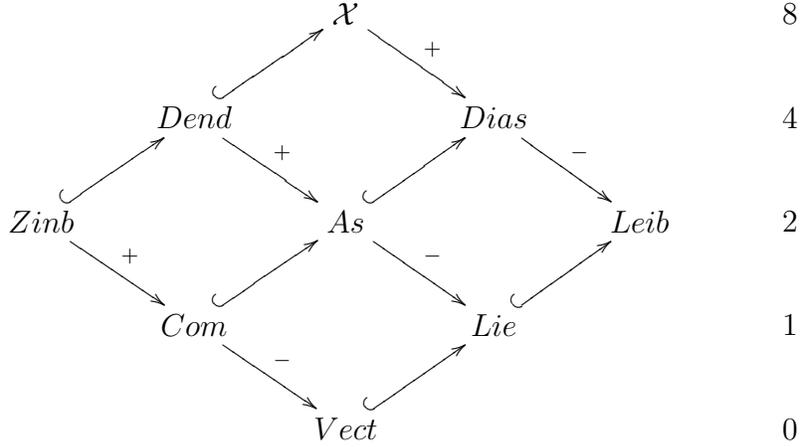
\[\xymatrix{
&& \XX\ar[rd]^{+} & & & \quad 8 \\
 &Dend\ar@{^{(}->}[ru]\ar[rd]^{+} & &Dias\ar[rd]^{-} & &\quad 4 \\
 Zinb\ar@{^{(}->}[ru]\ar[rd]^{+}& &As\ar@{^{(}->}[ru]\ar[rd]^{-} & &Leib &\quad 2\\
 &Com\ar@{^{(}->}[ru]\ar[rd]^{-}  & &Lie\ar@{^{(}->}[ru] & &\quad 1\\
 && Vect\ar@{^{(}->}[ru] && &\quad 0
}\]
\\

The numbers on the right side of the diagram indicate the dimension of the space of binary operations.

In other words we would like to find a notion of $\XX$-algebra whose operad is binary and quadratic, and satisfies the following properties:

\begin{enumerate}
\item the space of binary operations is 8 dimensional,
\item the operad is isomorphic to its dual (for Koszul duality),
\item a dendriform algebra is an $\XX$-algebra satisfying some symmetry,
\item any $\XX$-algebra gives, by some symmetrization of the operations, a diassociative algebra,
\item the functors deduced from the preceding two items make the upper square of the completed operadic butterfly commutative.
\end{enumerate}

The aim of this paper is to answer this question, and the answer is as follows: there are two solutions $\XX^+$ and $\XX^-$. 

An algebra of type $\XX^{\pm}$ has four generating operations denoted 
$ \begin{array}{cc}
\nwarrow  & \nearrow  \\
 \swarrow  & \searrow  \\
\end{array}
$
and $16=5\times 3 + 1 $ relations (we write $(\circ)\bullet$ instead of $(x\circ y)\bullet z$ ):

\[ \begin{array}{ccccccc}
(\nw) \nw = \nw (\nw) + \nw (\sw),   & (\sw) \nw = \sw (\nw), & (\nw) \sw + (\sw) \sw = \sw (\sw),  \\
(\nw) \nw = \nw (\se) + \nw (\ne),    & (\sw) \nw = \sw (\ne), & (\nw) \sw + (\sw) \sw = \sw (\se),     \\
(\ne) \nw = \ne (\nw) + \ne (\sw),    & (\se) \nw = \se (\nw), &(\ne) \sw + (\se) \sw = \se (\sw),   \\
(\nw) \ne = \ne (\ne) + \ne (\se),   & (\sw) \ne = \se (\ne), & (\nw) \se + (\sw) \se = \se (\se),  \\
(\ne) \ne = \ne (\ne) + \ne (\se),    & (\se) \ne = \se (\ne), & (\ne) \se + (\se) \se = \se (\se),  
\end{array}
\]
$$(\ne)\se - (\nw)\se = \pm( \ \nw(\sw) - \nw (\se)\ ) . \eqno{(16\pm)}$$
\\

A more conceptual way of describing  $\XX^{\pm}$ is as follows. Given two binary quadratic operads $\PP$ and $\QQ$ with a prefered basis for the space of generating operations, one can construct a new operad, denoted $\PP \sq \QQ$, whose set of generating operations  is the product of those of $\PP$ and $\QQ$ and the set of relations is also the product (in a certain sense) of those of $\PP$ and $\QQ$. This construction has already been used in particular cases  in \cite{A-L}, \cite{Leroux} (for ($\PP =\QQ$) and in \cite{E-G}. The operad $\XX^{\pm}$ is 
 $Dend \sq Dias$  quotiented by the relation $(16\pm)$.
\\

\noindent {\bf Content.} In the first part we recall the theory of Koszul duality for regular operads and we introduce the construction $\PP \sq \QQ$. In the second part we introduce the operads $Dend$ and $Dias$, and we compute $(Dend \sq Dias)^!$. In the third part we show that 
$$\XX^{\pm} := Dend \sq Dias / \textrm{ relation } (16\pm)$$
 completes the operadic butterfly. 


\noindent {\bf Convention.} All vector spaces are over the field $\mathbb K$. The tensor product over  
$\mathbb K$ of the two spaces $V$ and $W$ is denoted $V\t W$.

\section{Product of operads}\label{operads}

\subsection{Regular operads}\label{regular} In this paper we are dealing with algebras whose structure is defined by generating operations $(x,y)\mapsto x\circ_iy$ (with no symmetry) and relations of the form
$$\sum_{i,j} \alpha_{ij}\, (x\circ_iy) \circ _j z = \sum_{i,j} \beta_{ij}\, x\circ_i(y \circ _j z)\ ,\leqno (r)$$
where $ \alpha_{ij}$ and $ \beta_{ij}$ are scalars.
In these relations the variables stay in the same order (this is not the case for $Com, Lie, Zinb, Leib$). The associated operad, denoted $\PP$, is {\it binary} because the generating operations are binary. It is {\it  quadratic} because the relations involve monomials with two operations. It is {\it  regular} because the operations have no symmetry and, in the relations the variables stay in the same order. As a consequence of regularity the free $\PP$-algebra over the vector space $V$ is of the form 
$$\PP(V) = \oplus_{n\geq 1} \big( \PP_n \otimes {\mathbb K}[S_n]\big)\otimes_{S_n}V^{\t n}= \oplus_{n\geq 1} \PP_n \otimes V^{\t n},$$
 where $S_n$ is the symmetric group. The space $\PP_n$ is the space of (non-symmetric) $n$-ary operations.

Let us denote by $E=\PP_2$ the space of  (non-symmetric) binary operations. Let  $\{\circ_i\ \vert \  i\in I\}$ be a basis of $E$. The space of (non-symmetric) relations  $R$ is a subspace of $2E\t E$. The first summand corresponds to the parenthesizing $ (x\circ_iy) \circ _j z$ and the second summand to the parenthesizing $ x\circ_i(y \circ _j z)$. It will be helpful to denote by $ (\circ_i) \circ _j $ and $ \circ_i(\circ _j )$ the generating elements of $2E\t E$. Hence $R$ is the  subspace of $2E\t E$ generated by the vectors
$$r= \sum_{i,j} \alpha_{ij} (\circ_i) \circ _j - \sum_{i,j} \beta_{ij} \circ_i(\circ _j)\in 2E\t E\ $$
for each relation $(r)$.

One has $\PP_1 = {\mathbb K}$, since, up to multiplication of scalars, there is only one unary operation: the identity. One has $\PP_2 = E$  and 
$\PP_3 = 2E\t E/R$.

If an operad $\QQ$ is obtained from an operad $\PP$ by enlarging the space of relations with some more relations $(r)$, then the space of $n$-ary operations of $\QQ$ is a quotient of the space of $n$-ary operations of $\PP$. By abuse of language we will say that $\QQ$ is a quotient of $\PP$ by $(r)$ and we write $\QQ=\PP / (r)$ . 

The Koszul dual of the operad $\PP$, denoted $\PP^!$, is determined by $E^*:=\Hom(E, {\mathbb K})$ and $R^{\perp}$. Since we equipped $E$ with a basis, we can take the dual basis for $E^*$ and identify it with $E$. After this identification $R^{\perp}$ is described as follows: it is the space orthogonal to $R$ for the inner product $\langle -, - \rangle$ given on $E\t E \oplus E\t E$ by the matrix 
\begin{displaymath}
\left(\begin{array}{cc}
1 & 0 \\
0 & -1 
\end {array}\right)\ .
\end{displaymath}

In other words one has 
\begin{displaymath}
\left\{\begin{array}{cc}
\langle (\circ_i)\circ_j, (\circ_i)\circ_j \rangle = 1,& \\
\langle \circ_i(\circ_j), \circ_i(\circ_j) \rangle = -1,& \\
 \langle -, - \rangle=0 & \textrm{ in the other cases.}
\end {array}\right.
\end{displaymath}

Since the operations of $\PP$ and $\PP^!$ do not satisfy, in general, the same relations, it is necessary, sometimes, to distinguish between them. So we adopt the notation $\circ^*$ for the latter.

\begin{lemma}\label{orthogonality} Let $K$ be an index set and let $(r_k), k\in K$, be the relations defining $R$. An element
$$ \sum_{i,j} \alpha'_{ij} (\circ_i^*) \circ _j^* - \sum_{i,j} \beta'_{ij} \circ_i^*(\circ _j^*)\in 2E^*\t E^*\ $$
is in $R^{\perp}$ if and only if one has 
$$ \sum_{i,j} \alpha^{(k)}_{ij}\alpha'_{ij} =  \sum_{i,j} \beta^{(k)}_{ij}\beta'_{ij} $$
for all $k\in K$.
\end{lemma}
\noindent Proof.  This is an immediate translation of the definition of orthogonality. \hfill $\square$

\subsection{The square product of operads}\label{squareproduct} Let $\PP$ be a binary quadratic regular operad defined by binary operations denoted by $\circ_i$, and relations $(r)$ (cf.~\ref{regular}). Let $\QQ$ be another one with operations  $\bullet_k$ and relations $(r')$. We define the operad $\PP \sq \QQ $ by the operations $\circ_i \bullet_k$ (product of the two sets of operations), and relations $(r,r')$ given by
$$\sum_{i,j,k,l}\alpha_{ij}\alpha'_{kl} (\circ_i \bullet_k) \circ_j\bullet_l = 
\sum_{i,j,k,l}\beta_{ij}\beta'_{kl} \circ_i \bullet_k ( \circ_j\bullet_l)  .$$
So, if $\PP$ is defined by $m$ relations and $\QQ$ by $m'$ relations, then $\PP \sq \QQ $ is defined by $mm'$ relations.

It is immediate to see that the construction $\sq$ is associative, commutative, and its neutral element is the operad $As$. Indeed, $As$ has only one operation $\cdot$ satisfying $(\cdot) \cdot = \cdot (\cdot)$ .

\begin{prop}\label{squaredual} Let $\PP$ and $\QQ$ be binary quadratic regular operads. The operad 
$ (\PP \sq \QQ)^! $ is a quotient of the operad $\PP^! \sq \QQ^!$, so there is a natural forgetful functor of categories of algebras:
$$ (\PP \sq \QQ)^! {\rm -alg} \longrightarrow (\PP^! \sq \QQ^!){\rm -alg} .$$
\end{prop}

\noindent Proof. The set of generating operations is the same in both cases, only the space of relations is different. For $ (\PP \sq \QQ)^! $ it is $T^{\perp}$, where $T$ is the space generated by the relations $(r,r')$ (cf. \ref{squareproduct}). For $\PP^! \sq \QQ^!$ it is the space $S$ generated by the relations $(s,s')$, where $s$ is orthogonal to all the relations $r$, and $s'$ is orthogonal to all the relations $r'$.

In order to complete the proof it is sufficient to prove that $S$ is included into $T^{\perp}$. Indeed the expected functor would then simply be the forgetful functor. The space $S$ is included into $T^{\perp}$ if and only if $\langle S, T \rangle = 0$. Let us check this equality for the relations $(s,s')$ and $(r,r')$.

We denote by $\alpha, \beta$ the structure constants of $r$ as in \ref{regular} and by  $\gamma, \delta$ the structure constants of $s$  (same thing with a prime for $r'$ and $s'$). By Lemma \ref{orthogonality}, it suffices to prove the equality
$$\sum_{i,j,k,l}\alpha_{ij}\alpha'_{kl}\gamma_{ij}\gamma'_{kl}=   \sum_{i,j,k,l}\beta_{ij}\beta'_{kl} \delta_{ij}\delta'_{kl}  .$$
The lefthand summand is equal to 
$$(\sum_{i,j}\alpha_{ij}\gamma_{ij}) (\sum_{k,l} \alpha'_{kl}\gamma'_{kl})$$
and similarly for the righthand side. Since the relations $r$ and $s$ (resp. $r'$ and $s'$) are orthogonal one has
$$\sum_{i,j}\alpha_{ij}\gamma_{ij}= 
\sum_{i,j}\beta_{ij}\delta_{ij}$$
and similarly with a prime. This proves the equality and so $\langle S, T \rangle = 0$.\hfill $\square$

\section{Dendriform and diassociative algebras}\label{dend-dias} We recall from \cite{Lo} the definition of dendriform algebra (with a slight change of notation) and of diassociative algebra.

\subsection{Definition}\label{dend} A \emph {dendriform} algebra $A$ is a vector space equipped with two operations denoted $\wedge$ and $\vee$ satisfying the following relations:
\begin{displaymath}
\left\{\begin{array}{rclr}
(x\wedge y)\wedge z &=& x\wedge (y\wedge z) + x\wedge (y\vee z)  ,&(i) \\
(x\vee y)\wedge z &=& x\vee (y\wedge z) ,&(ii)  \\
(x\wedge y)\vee z  + (x\vee y)\vee z &=& x\vee (y\vee z).&(iii) 
\end{array}
\right.
\end{displaymath}

\subsection{Definition}\label{dend} A \emph{diassociative} algebra (or associative dialgebra) $A$ is a vector space equipped with two operations denoted $\g$ and $\d$ satisfying the following relations:
\begin{displaymath}
\left\{\begin{array}{rclcr}
(x\g y)\g z &=& x\g (y\g z) ,&\qquad&(1)  \\
(x\g y)\g z &=& x\g (y\d z) , &\qquad&(2) \\
(x\d y)\g z &=& x\d (y\g z) , &\qquad&(3) \\
(x\g y)\d z &=& x\d (y\d z) , &\qquad&(4) \\
(x\d y)\d z &=& x\d (y\d z) . &\qquad&(5) 
\end{array}
\right.
\end{displaymath}
It was shown in \cite{Lo} that the associated operads are dual to each other via the identification $\wedge^* = \g $ and $\vee^* = \d$. So we have $Dend^!=Dias$ and $Dias^!=Dend$. In fact the reader can check it immediately from the description of duality for regular operads given in Lemma \ref{orthogonality}.

\subsection{The Operad $Dend\sq Dias$ and its dual}\label{dualdenddias}
Consider the operad $Dend\sq Dias$. We denote its set of generating operations by
$$ \nw:=(\wedge, \g),\  \ne:=(\wedge, \d),\   \se:=(\vee, \d),\  \sw:=(\vee, \g)\ .$$
Since $Dend$ has 3 relations and since $Dias$ has 5 relations, $Dend\sq Dias$ has 15 relations which read as follows:

\[ \begin{array}{ccccccc}
(\nw) \nw = \nw (\nw) + \nw (\sw),   & (\sw) \nw = \sw (\nw), & (\nw) \sw + (\sw) \sw = \sw (\sw),  \\
(\nw) \nw = \nw (\se) + \nw (\ne),    & (\sw) \nw = \sw (\ne), & (\nw) \sw + (\sw) \sw = \sw (\se),     \\
(\ne) \nw = \ne (\nw) + \ne (\sw),    & (\se) \nw = \se (\nw), &(\ne) \sw + (\se) \sw = \se (\sw),   \\
(\nw) \ne = \ne (\ne) + \ne (\se),   & (\sw) \ne = \se (\ne), & (\nw) \se + (\sw) \se = \se (\se),  \\
(\ne) \ne = \ne (\ne) + \ne (\se),    & (\se) \ne = \se (\ne), & (\ne) \se + (\se) \se = \se (\se). 
\end{array}
\]

\centerline {The tableau}

\begin{prop}\label{prop:denddiasdual} The operad of $(Dend\sq Dias)^!$-algebras is isomorphic to the quotient of the operad $Dend\sq Dias$ by the following two relations
\begin{displaymath}
\left\{\begin{array}{c}
(\ne)\se - (\nw)\se =0\ ,\\
0=\nw(\sw) - \nw(\se)\ .
\end{array}
\right.
\end{displaymath}
\end{prop}

\noindent Proof. The 15 relations of $Dend\sq Dias$ are linearly independent, hence the space of relations is of dimension 15. Its orthogonal is of dimension $2\times 4^2 - 15 = 17$. Since $Dend^!=Dias$ and $Dias^!=Dend$ one has $Dend^!\sq Dias^! = Dias \sq Dend \cong Dend \sq Dias$. By Proposition \ref{squaredual} $(Dend\sq Dias)^!$ has, at least,  the 15 relations of $ Dend \sq Dias$ (up to isomorphism) as relations plus two more (linearly independent). Let us show that the two extra operations are the ones indicated.

We first make explicit the isomorphism:
\begin{displaymath}
\begin{array}{c}
\nw^* =(\wedge,\g)^* = (\wedge^*,\g^*) = (\g, \wedge) \cong \nw\ ,\\
\ne^* =(\wedge,\d)^* = (\wedge^*,\d^*) = (\g, \vee) \cong \sw\ ,\\
\se^* =(\vee,\d)^* = (\vee^*,\d^*) = (\d, \vee) \cong \se\ ,\\
\sw^* =(\vee,\g)^* = (\vee^*,\g^*) = (\d, \wedge) \cong \ne\ .\\
\end{array}
\end{displaymath}
So, in order to write down the 15 relations of $(Dend\sq Dias)^!$ with this basis of generating operations we have to exchange the operations $\ne$ and $\sw$ in the tableau of the 15 relations. Let us call it the dual tableau. Under this isomorphism the two extra relations read:
\begin{displaymath}
\left\{\begin{array}{c}
(\sw)\se - (\nw)\se =0\ ,\\
0=\nw(\ne) - \nw(\se)\ .
\end{array}
\right.
\end{displaymath}
Let us check that these two relations are orthogonal with the 15 relations of $Dend\sq Dias$.
In 28 cases the verification is immediate because the involved operations are all different. The remaining two cases 
are
\begin{displaymath}
\begin{array}{c}
\langle\ (\sw)\se - (\nw)\se\ , \ (\nw)\se + (\sw)\se - \se(\se)\ \rangle = +1-1=0 ,\\
\langle\ \nw(\ne) - \nw(\se)\ , \ -(\nw)\nw + \nw(\se) + \nw(\ne)\ \rangle = -1+1=0 .
\end{array}
\end{displaymath}
It is straightforward to check that these 17 relations are linearly independent, therefore we have a complete presentation of $(Dend\sq Dias)^!$. Under the inverse isomorphism of the one described above we get the expected result.\hfill $\square$

\section {The missing operad}\label{missing}
In this section we complete the operadic butterfly by constructing the operad $\XX$.

\begin{thm}\label{ex} Let $\XX^+$ (resp. $\XX^-$) be the operad $Dend\sq Dias$ quotiented by the relation (16+) (resp. (16-)):
\begin{displaymath}
\begin{array}{rcl}
(\ne)\se - (\nw)\se &=& +\nw(\sw) - \nw(\se)\ ,\qquad (16+)\\
(\ne)\se - (\nw)\se &=& -\nw(\sw) + \nw(\se)\ .\qquad (16-)
\end{array}
\end{displaymath}
The operad $\XX^{\pm}$ completes the operadic butterfly, that is, it satisfies the following properties:
 
\begin{enumerate}
\item the space of binary operations is 8 dimensional,
\item the operad is isomorphic to its dual (for Koszul duality),
\item a dendriform algebra is an $\XX$-algebra satisfying some symmetry,
\item any $\XX$-algebra gives, by some symmetrization of the operations, a diassociative algebra,
\item the functors deduced from the preceding two items make the upper square of the completed operadic butterfly commutative.
\end{enumerate}
\end{thm}

\noindent Proof. Let us put $r:= (\ne)\se - (\nw)\se $ and $s:= \nw(\sw) - \nw(\se)$. We are looking for scalars $\alpha$ and $\beta$
such that the operad $Dend\sq Dias/ \alpha r + \beta s$ is self-dual. The dual of $r$ (resp. $s$) is 
$r^*= (\sw)\se - (\nw)\se$ (resp. $s^* = \nw(\ne) - \nw(\se)$ ). To get self-duality of the quotient operad we need to have 
$-\alpha^2 + \beta ^2 = 0$. Hence $\alpha$ and $\beta$ have to be equal up to sign. So the sixteenth relation is either $r=s$, called $(16+)$ or  $r=-s$, called $(16-)$. It is immediate to check that in both cases the 16 ($=3\times 5 + 1)$) relations are linearly independent. Hence we have proved (1) and (2). It turns out that both solutions satisfy also the other requirements as we prove now.
\\

\noindent (3) Let $(A, \wedge, \vee)$ be a dendriform algebra. Define $\nw := \wedge=: \ne$ and $\sw := \vee=: \se$. Let us prove that $(A, \nw, \ne, \sw, \se)$ is an $\XX^{\pm}$-algebra. The first fifteenth relations are fullfilled since relation $(n,a)$ is a consequence of relation $(a)$ for $n=1,\ldots , 5$ and $a = i, ii, iii$. The relation $(16\pm)$ is also fullfilled since both sides of the equality are $0$: the left side because $\nw =\ne$ and the right side because $\sw =\se$.
\\

\noindent (4) Let $(A, \nw, \ne, \sw, \se)$ be an $\XX^{\pm}$-algebra. Define 
$$x\g y := x\nw y + x\sw y \quad   {\rm and }\quad 
 x\d y := x\ne y + x\se y.$$
 Adding relations $(n,i), (n,ii)$ and $(n,iii)$ shows that relation $(n)$ is fullfilled for 
$n=1,\ldots , 5$.
\\

\noindent (5) Starting with a dendriform algebra $(A, \wedge, \vee)$, we get a diassociative algebra  $(A, \g, \d)$ such that 
$$x\g y = x\nw y + x\sw y =  x\wedge y + x\vee y = x\ne y + x\se y = x \d y\ .$$
Therefore the composite $Dend \hookrightarrow \XX^{\pm}  \stackrel{+}{\to} Dias$ is equal to the composite
$Dend  \stackrel{+}{\to}  As \hookrightarrow Dias$ as expected. \hfill $\square$
\\

\noindent {\bf Remark.} In fact our proof shows that these two solutions are the only quotients of $Dend\sq Dias$ which satisfy all the required properties.
\\

\subsection{Question}  The operads $Com, Lie, As, Zinb, Leib, Dend, Dias$ are Koszul operads, that is the associated Koszul complex is acyclic (cf. \cite {G-K}, \cite{L-P}, \cite{Lo}). Is also $ \XX^{+} $(resp.  $ \XX^{-} $) a Koszul operad ?

Since the operad $\XX$ ($=\XX^+$ or $\XX^-$) is regular and self-dual, the criterion to ensure Koszul duality takes the following form. Let $\XX_n$ be the homogeneous part of degree $n$ in the free $\XX$-algebra on one generator. Then for each integer $k\geq 1$ there is a finite chain complex

\[\begin{array}{c}
0\to \XX_k\t\XX_1\t \cdots \t \XX_1 \to \cdots \hfill \\
\qquad \qquad \to  \bigoplus_{m_1+ \cdots +m_n=k}\XX_n\t\XX_{m_1}\t \cdots\t \XX_{m_n}\to\cdots \qquad \qquad\\
\hfill  \cdots \to \XX_1\t\XX_k\to 0.
\end{array}
\]

The acyclicity of these complexes for $k>1$ imply Koszulity. It would imply that $\XX_n$ is of dimension $4^{n-1}$.

\end{document}